\documentclass{amsart}
\usepackage{amssymb}

\newtheorem{thm}{Theorem}[section]

\newtheorem{lemma}[thm]{Lemma}
\newtheorem{prop}[thm]{Proposition}
\theoremstyle{definition}
\newtheorem{defn}[thm]{Definition}
\theoremstyle{remark}
\newtheorem{rem}[thm]{Remark}

\numberwithin{equation}{section}
\newtheorem*{hyp}{Hypothesis}

\newcommand{\GS}{Gr\"obner-Shirshov}
\newcommand{\vb}{\mathord{\hskip 1pt|\hskip 1pt}}

\newcommand{\CD}{Composition-Diamond}
\begin{document}

\title{Gr\"obner-Shirshov basis for the braid semigroup  }

\author[Bokut]{L. A. Bokut$^*$}
\address{Sobolev Institute of Mathematics, Novosibirsk 630090, Russia; =
and
Chang Jung Christian University, Kway Jen, Tainan 711, Taiwan}
\email{bokut@math.nsc.ru}
\thanks{$^*$Supported in part by the Russia's Fund of
    Fundamental Research.}

\author[Fong]{Y. Fong}
\address{Department of Mathematics, National Cheng-Kung University,
Tainan , Taiwan} \email{yfong@mail.ncku.edu.tw}

\author[Ke]{W.-F. Ke}
\address{Department of Mathematics, National Cheng-Kung University,
Tainan , Taiwan} \email{wfke@mail.ncku.edu.tw}

\author[Shiao]{L.-S. Shiao}
\address{Chang Jung Christian University, Kway Jen, Tainan 711,
Taiwan} \email{lsshiao@mail.cju.edu.tw}

 \begin{abstract}
 We found Gr\"obner-Shirshov basis for the braid semigroup $B^+_{n+1}$.
 It gives a new algorithm for the solution of the word problem for the braid
 semigroup and so for the braid group.
  \end{abstract}

\maketitle

\section{Introduction}

In this paper, we are continuing to use a language of the \GS\
bases for some semigroups presented by defining relations (in
papers \cite{BoS01}, \cite{BoSh01} some groups have been treated
as semigroups). Traditionally in this situation, one use some other
languages like  elementary transformations language
(\cite{Ad66},\cite{AD00}), Newman's Diamond lemma language
(\cite{Ne42}), the rewriting systems language (\cite{BoOt93} ), or
some else (\cite{KhS95}). For example, P.M. Cohn used Diamond
lemma for semigroups as early as in 1956 (\cite{Co56a},
\cite{Co56b}). After him first author used the same language to
prove some embedding theorems for semigroups
\cite{Bo63a},\cite{Bo63b}. Also first author used a language of
rewriting systems in order to prove that multiplicative semigroups
of some rings is  embeddable into groups \cite{Bo69a},
\cite{Bo69b} (see also  \cite{Bo80}, \cite{Bo96}).
 But here we formulate our main results in terms of \GS\ bases.

During the last few years the \GS\ bases method has been used to
study some important examples of associative algebras (quantum
enveloping algebras of the type $A_n$, Hecke algebras of type
$A_n$), Lie algebras (simple Lie algebras of the types $A_n, B_n,
C_n, D_n,$ $ G_2, F_4, E_6, E_7, E_8$, Kac-Moody algebras of the
types $A_n^{(1)},B_n^{(1)},C_n^{(1)},D_n^{(1)}$), Lie superalgebras
(simple Lie superalgebras of
the types $A_n, B_n, C_n, D_n$), modules (irreducible modules over
simple Lie algebras of the type $A_n$, Specht modules over Hecke
algebras of the type $A_n$), and groups (\GS\ bases for Coxeter
groups of types $A_n, B_n, D_n$) \cite{BoK96} -- \cite{BoS01},
\cite{KaLe00}--\cite{KangPreprint}, \cite{Po99a} -- \cite{Po99d}. Also this method
was invented for associative conformal algebras \cite{BoFK1},
\cite{BoFK2}.

In recent paper \cite{BoSh01},  two of us  applied the \GS\ bases method
to study some groups by P.S. Novikov \cite{No55} and W.W. Boone
\cite{Boone59} with the algorithmically unsolvable word problem.
They revised some papers \cite{Bo66}, \cite{Bo67a}, \cite{Bo84},
\cite{Bo85} (see also \cite{Bo86}, \cite{BoK94}) that used a
rewriting systems language (between lines  of HNN-extension
theory) using the very language of \GS\ bases. It means that they
eliminated any group theory in proving Novikov's and Boone's
theorems on unsolvability.

Some other activities on \GS\ (non-commutative Gr\"obner,
standard) bases  may be found in papers by T. Mora \cite{Mo94},
V.N. Latyshev \cite{La00}, V. Ufnarovski \cite{Uf98}, A.A.
Mikhalev and A.A. Zolotyh \cite{MiZ98}, A.A. Mikhalev and E.
Vasilieva \cite{MiV00}, V.P. Gerdt and V.V. Kornyak
\cite{GerdtK96}.

A great  interest has been done to the word and conjugacy problems
for  braid groups $B_{n+1}$. Let us mention a pioneering paper by
E. Artin \cite{Ar26}, \cite{Ar47}, papers by A.A. Markov
\cite{Ma45}, A.F. Garside \cite{Ga69}, P. Deligne \cite{Del72}, W. Thurston \cite{Th88},
K. Tatsuoka \cite{Ta93},
E.A. Elrifai and H.R. Morton
\cite{EM94}, P. Dehornoy \cite{De97}, J. Birman, K.H. Ko, and S.J.
Lee \cite{BKL98}. In those papers , a fundamental role of a semigroup of positive words
of $B_{n+1}$ is proved to be crucial. Actually this semigroup
depends of generators of $B_{n+1}$. Let us denote this semigroup
in Artin generators $a_i,1\leq i\leq n$ by $B_{n+1}^+$ and call
the braid semigroup. It has the following
defining relations as a semigroup:
\begin{align*}
&a_{i+1}a_ia_{i+1}=a_ia_{i+1}a_i,\\
&a_ia_j=a_ja_i, i-j\geq 2.
\end{align*}

In this paper, we found \GS\ basis of $B_{n+1}^+$. It gives
a new solution of the word problem for the braid semigroup and as the
result for the braid group. We also
 study a more general concept of braid semigroups. The braid
semigroup $B_{n+1}^+$ is  a braid semigroup of type
$A_n$.
 We formulate a hypothesis on \GS\ bases
of any braid semigroup. It is a modification of a hypothesis on
\GS\ bases of Coxeter groups \cite{BoS01}. The previous example
support our hypothesis.

Let us mentioned that in Birman-Ko-Lee generators $a_{ts}, 1\leq
t<s\leq n+1$ a semigroup of positive braids of $B_{n+1}$
have the following defining relations \cite{BKL98}:
\begin{align*}
&a_{ts}a_{rq}=a_{rq}a_{ts}, (t-r)(t-q)(s-r)(s-q)>0,\\
&a_{ts}a_{sr}=a_{tr}a_{ts}=a_{sr}a_{tr} .
\end{align*}

It would be interesting to fined a \GS\ basis of this semigroup
too. It may have some connections with a new Public-Key
Cryptosystem \cite{Ko00}.

\section{\GS\ bases }

Let $\mathbf{X}$ be a linearly ordered set,  $k$ be a field,
$k\langle \mathbf{X}\rangle$ be the free associative algebra over
$\mathbf{X}$ and $k$. On the set $\mathbf{X}^*$ of words we impose
a monomial well order $>$ (i.e. a well order that  agrees with the
concatenation of words). Any polynomial $f\in k\langle
\mathbf{X}\rangle$ has the leading word $\bar f$. We say that $f$
is monic if $\bar f$ occurs in $f$ with coefficient $1$. By a
composition of intersection $(f,g)_w$ of two monic polynomials
relative to some word $w$, such that $w = \bar fb= a\bar g$,
$\deg( \bar f) + \deg(\bar g) > \deg (w)$, one means the following
polynomial
$$
                            (f,g)_w = fb - ag.
$$
By composition of including $(f,g)_w$  two monic polynomials,
where $w = \bar f = a\bar gb$, one means the following polynomial
$$
                           (f,g)_w = f - agb.
$$
In the last case the transformation
$$
                     f\mapsto (f,g)_w = f -agb
$$
is called the elimination of the leading word (ELW) of g in f.

In both cases the word $w$ is called an $ambiguity$ of polynomials
(relations) $f,g.$

A composition $(f,g)_w$ is called trivial relative to some
$\mathbf{S} \subset k\langle \mathbf{X}\rangle$ ($(f, g)_w\equiv
0 $ $ mod (\mathbf{S}, w))$, if
$$
            (f,g)_w = \sum \alpha_i a_it_i b_i,
$$
where   $t_i \in \mathbf{S}$, $a_i, b_i \in  \mathbf{X}  ^*$, and
$\bar {a_it_ib_i} < w$. In particular, if $(f,g)_w$ goes to zero
by the ELW's of $\mathbf{S}$ then $(f,g)_w$ is trivial relative to
$\mathbf{S}$.

We assume for some polynomials $f_1$ and $f_2$ that
$$
               f_1 \equiv f_2 mod (\mathbf{S}, w),
$$
if
$$
f_1-f_2 \equiv 0 mod (\mathbf{S},w).
$$

A subset $\mathbf{S}$ of $k\langle \mathbf{X}\rangle$ is called a
\GS\ basis if any composition of polynomials from $\mathbf{S}$ is
trivial relative to $\mathbf{S}$.

By $\langle \mathbf{X}\vb \mathbf{S}\rangle$, the algebra with
generators $\mathbf{X}$ and defining relations $\mathbf{S}$, we
will mean the factor-algebra of $k\langle \mathbf{X}\rangle$ by
the ideal generated by $\mathbf{S}$.

The following lemma goes back to the Poincare - Birkhoff-Witt
(PBW) theorem, the Diamond Lemma of M.H.A. Newman \cite{Ne42}, the
Composition Lemma of A.I. Shirshov (\cite{Sh62}) (see also
\cite{Bo72}, \cite{Bo76}, where this Composition Lemma was
formulated explicitly and in a current form for Lie algebras and
associative algebras respectively), the Buchberger's Theorem
(\cite{Bu65}) (published in \cite{Bu70}), the Diamond Lemma of G.
Bergman \cite{Be78} (this lemma was also known to P.M. Cohn, see,
for example, \cite{Co66}, \cite{Co81} and some historical comments
to Chapter ``Gr\"obner bases" in \cite{Ei95}):

\begin{lemma}(Composition-Diamond Lemma)

   $\mathbf{S}$ is a \GS\ basis if and
only if the set
$$
       \{u \in \mathbf{X}^* \mid u \neq a\bar fb, \mbox{ for any }f\in
\mathbf{S}\}
$$
of $\mathbf{S}$-reduced words consists of  a linear basis of the
algebra $\langle \mathbf{X}\vb \mathbf{S}\rangle$.
\end{lemma}

While $\mathbf{S}$ is fixed, we will call $\mathbf{S}$-reduced
words just reduced words. Actually we would like to introduce a
new term for this basis - the $ PBW-basis$ for the algebra relative
to the \GS\ basis of relations of this algebra. It would reflect to
the fact that as we have mentioned the Composition-Diamond lemma
goes back to the Poincare-Birkhoff-Witt (PBW) theorem as well.
In fact the PBW theorem is a particular and the most important up
to now case of \CD\ lemma when $\mathbf{S}= \{x_ix_j-x_jx_i-\sum
x_{ij}^k, i>j\}$ and $[x_ix_j]=\sum x_{ij}^k, i>j$ is the
multiplicative table of some Lie algebra. In this case
$PBW(\mathbf{S})=\{x_{i_1}x_{i_2}\dots x_{i_k}, i_1\leq i_2\dots
\leq i_k \}.$

If a subset $\mathbf{S}$ of $k\langle \mathbf{X}\rangle$ is not a
\GS\ basis then one can add to $\mathbf{S}$ all nontrivial
compositions of polynomials of
 $\mathbf{S}$, and continue this process (infinitely) many times  in
order to have a
\GS\  basis $\mathbf{S}^{\rm comp}$ that contains $S$. This
procedure is called the Buchberger - Shirshov algorithm
\cite{Sh62}, \cite{Bu65}, \cite{Bu70}.

A \GS\ basis $\mathbf{S}$ is called $minimal$  if any
$s\in \mathbf{S}$ is a linear combination of
$\mathbf{S}\setminus\{s\}$ - reduced words. Any ideal of $k\langle
\mathbf{X} \rangle$ with a fix monomial order has the unique minimal \GS\ basis.

If $\mathbf{S}$ is a set of ``semigroup relations" (that is,
polynomials of the form $u - v$, where $u,v \in \mathbf{X}^*$),
then any nontrivial composition will have the same form. As a
result the set $\mathbf{S}^{\rm comp}$ also consists of semigroup
relations.

Let $A = smg\langle \mathbf{X}\vb \mathbf{S}\rangle$  be a
semigroup presentation. Then $\mathbf{S}$ is a subset of $k\langle
\mathbf{S}\rangle$ and one can find a \GS\ basis $ \mathbf{S}^{\rm
comp}$. The last set does not depend on $k$, and consists of
semigroup relations. We will call $\mathbf{S}^{\rm comp}$ a \GS\
basis of $A$. It is the same as a \GS\ basis of the semigroup
algebra $kA = \langle \mathbf{X}\vb\mathbf{S}\rangle$.

 Let $\{u_i\}$ be a set of words in an alphabet. By $V(u_i)$ we
understand any
word in $\{u_i\}$. Let $\{w_i\}$ be some other set of words with a
surjection $u_i \mapsto w_i$. Then by $V(w_i)$ we understand the
same word  $V$, but in the set $\{w_i\}$.

\section{A hypothesis}

Here we will formulate a general hypothesis  on \GS\ bases for any
braid semigroup.

We will present Coxeter group $W$  relative to Coxeter matrix $M$
in the following way. Let $S$ be a linear ordered $l$-elements
set, $M = (m_{ss'})$  be $l \times l$ Coxeter matrix. Then
$$
   W = smg\langle S \mid s^2 = 1, (ss' )^{m_{ss'}} = 1, s \neq s',
\forall s, s' \in S \mbox{ and finite }m_{ss'}\rangle.
$$
Let us define  for finite $m_{ss'}$,
\begin{align*}
m(s, s')  &= ss' \dots &&\mbox{(there are $m_{ss'}$ alternative
    letters $s$, $s'$)},\\
(m-1) (s, s') &= ss' \dots &&\mbox{(there are $m_{ss'}-1$
    alternative letters $s$, $s'$)},
\end{align*}
and so on.  For example, if $m_{ss'} = 2$, then $m(s, s') = ss'$,
$(m-1)(s, s') = s$. If $m_{ss'} = 3$, then $m(s, s') = ss's$,
$(m-1)(s, s') = ss'$.

In this notations defining relations of W can be presented in the
form
\begin{equation}\label{E1}
                 s^2 = 1,\,  m(s, s') = m(s', s),\,  s > s',
\end{equation}
for all $s, s' \in S$ and finite $m_{ss'}$.

\begin{defn}
The following semigroup will be called the braid semigroup $W^+$
relative to the Coxeter matrix $M$:
\begin{equation}\label{E2}
W^+= smg\langle s_1, \dots, s_l\vb  m(s, s') = m(s', s),\,  s >
s', s,s'\in S, m_{ss'}<\infty \rangle.
\end{equation}

\end{defn}

Also we will use the following notation:
$$
(m,i)(s,s')
$$
is a result of removing first $i$ letters in $m(s,s')$, where
$1\leq i\leq m.$

Let us call two words in $S$ equivalent if they are equal module
commutativity relations from (\ref{E1}).  To be more prices, it
means that they are equal in the so called free partially
commutative semigroup (algebra) generated by $S$ and with
commutativity relations from (\ref{E1}) \cite{BoS01}.  Two
relations $a = b$ and $c = d$ of $W$ are called equivalent if $a$,
$b$ are equivalent to $c$, $d$ respectively.

\begin{hyp}
 A \GS\ basis $\mathbf{T}$ of $W^+$ include of initial relations
(\ref{E2}) and relations that are equivalent to  the following
ones:
\begin{equation}\label{E3}
\begin{split}
 &(m-i_0)(s, s')(m-i_1)(s_1, s_2) \dots (m-i_k)(s_{2k-1},
s_{2k})m(s_{2k+1},
  s_{2k+2})\\
 &\quad = m(s', s)(m,i_0)(s_2, s_1) \dots (m,i_{k-1})(s_{2k},
s_{2k-1})(m,i_k)
(s_{2k+2}, s_{2k+1}),
\end{split}
\end{equation}
where $s > s'$, $s_1 < s_2$, \dots, $s_{2k-1} <  s_{2k}$,
$s_{2k+1} < s_{2k+2}$,  and
\begin{align*}
&(m-i_k)(s_{2k-1}, s_{2k})m(s_{2k+2},s_{2k+1})\doteq
m(s_{2k-1}, s_{2k})(m,i_k)(s_{2k+2},s_{2k+1}),\\
 &(m-i_{k-1})(s_{2k-3},
s_{2k-2})m(s_{2k},s_{2k-1})\doteq m(s_{2k-3},
s_{2k-2})(m,i_{k-1})(s_{2k},s_{2k-1}),\\
 &.\\
 &. \\
 &. \\
&(m-i_1)(s_1, s_2)m(s_4,s_3)\doteq
m(s_1, s_2)(m,i_1)(s_4,s_3),\\
 &(m-i_0)(s,
s')m(s_2,s_1)\doteq m(s, s')(m,i_0)(s_2,s_1).
\end{align*}

Above $\doteq$ means the equality in the free semigroup.

Let us fix left hand and rite hand sides of (\ref{E3}). If $A=B$
has a form (\ref{E3}) then a transformation
$$
A\rightarrow B
$$
will be  called a positive chain, while $B\rightarrow A$ will be
called negative.

Any other relation of $\mathbf{T}$ is equivalent to the following
one:
\begin{equation}\label{E4}
X=Y,
\end{equation}
where $X$ goes to some $X'$ by a series of negative chains,
$$
X\rightarrow X_1\rightarrow \dots \rightarrow X_k=X',
$$
(in the free  partially commutative semigroup), $X'=AY_1$, and
$Y=BY_1$, where $A\rightarrow B$ is a positive chain.
\end{hyp}

\begin{rem}
Let us compere this hypothesis with our hypothesis on \GS\ bases
of Coxeter groups \cite{BoS01}. Relations (6.2) in \cite{BoS01} is
just relations (\ref{E3}) with $i_0=i_1= \dots =i_k=1$ and with an
extra condition that all neighbor pairs $(s,s'), (s_1,s_2), \dots
,(s_{2k+1},s_{2k+2})$ are different. In particular we don't expect
negative chains for Coxeter groups at all.
\end{rem}

\section{Braid semigroup of the type $A_n$}

Let $B_3^+$ is the braid semigroup with 2 generators
$a_1,a_2, a_2a_1a_2=a_1a_2a_1.$ We will assume that $a_1<a_2$, and
$u<v$ for words $u,v$ will mean the $deg-lex$ order (first we
compere words by the degree (length) and then lexicographically).

\begin{lemma}\label{L1}
The minimal \GS\ basis of $B_3^+$ consists of the initial
relation and relations
\begin{equation}\label{E5}
a_2a_1^la_2a_1=a_1a_2a_1^2a_2^{l-1}, l\geq 2
\end{equation}
of a form (\ref{E3}).
\end{lemma}

\begin{proof}
First of all (\ref{E5}) is follows from Garside's formulas
\cite{Ga69}, Lemma 2:
$$
a_1\Delta_2=\Delta_2a_2, a_2\Delta_2=\Delta_2a_1,
$$
where $\Delta_2=a_1a_2a_1.$

Secondly (\ref{E5}) has a form (\ref{E3}):
$$
(m-1)(a_2,a_1)(m-2)(a_1,a_2)(m-2)(a_1,a_2) \dots
(m-2)(a_1,a_2)m(a_1,a_2)=
$$
$$
m(a_1,a_2)(m,1)(a_2,a_1)(m,2)(a_2,a_1)\dots
(m,2)(a_2,a_1)(m,2)(a_2,a_1).
$$

Now we need to proof that any composition of relations from the
lemma is trivial. There are some compositions of intersections
relative to the following   ambiguities $w$:
\begin{align*}
&a_2a_1a_2a_1a_2, a_2a_1^la_2a_1a_2, l\geq
2, a_2a_1a_2a_1^la_2a_1, l\geq 2,\\
&a_2a_1^la_2a_1^ka_2a_1, l,k\geq 2.
\end{align*}

Let us check the last composition for an example:
\begin{align*}
&(f,g)_w=(a_2a_1^la_2a_1-a_1a_2a_1^2a_2^{l-1})a_1^{k-1}a_2a_1-\\
&a_2a_1^l(a_2a_1^ka_2a_1-a_1a_2a_1^2a_2^{k-1})=\\
&a_2a_1^{l+1}a_2a_1^2a_2^{k-1}-a_1a_2a_1^2a_2^{l-1}a_1^{k-1}a_2a_1\equiv\\
&a_1a_2a_1^2a_2^la_1a_2^{k-1}-a_1a_2a_1^2a_2^{l-1}a_1^{k-1}a_2a_1\equiv
0.
\end{align*}

All other compositions can be figure out in the same way.

\end{proof}

Let us introduce the following notations that are  analogous to
that of \cite{BoS01}:
\begin{align*}
&a_{ij}=a_ia_{i-1} \dots a_j,i\geq j, a_{ii}=a_i,a_{ii+1}=1,\\
& a_{ij}(\mathbf{p})=a_i^{p_0}a_{i-1}^{p_1} \dots
a_j^{p_{i-j}},i\geq j,a_{ii}( \mathbf{p})=a_i^{p_0},a_{ii+1}(
\mathbf{p})=1.
\end{align*}
 where
$\mathbf{p}=(p_0,\dots, p_{i-j})$, all $p'$s are $\geq 1$. We will
write $\mathbf{p}> \mathbf{1}$, if at least one of $p$'s is
greater that 1. Otherwise we will write $\mathbf{p}= \mathbf{1}.$

 It is easy to see the
following equality in $B_n^+$:
\begin{equation}\label{E6}
\begin{split}
&a_{ij-1}a_{i+1j-1}=a_{i+1j-1}a_{i+1j},j\leq i+1.\\
&a_{ij-1}( \mathbf{p})a_{i+1j-1}=a_{i+1j-1}a_{i+1j}( \mathbf{p}),
j\leq i+1.
\end{split}
\end{equation}
In fact the first of this equalities has a form:
\begin{equation}\label{E7}
\begin{split}
&(m-1)(a_i,a_{i+1})(m-1)(a_{i-1},a_i) \dots
(m-1)(a_j,a_{j+1})m(a_{j-1},a_j)=\\
&m(a_{i+1},a_i)(m,1)(a_{i-1},a_i)\dots
(m,1)(a_{j+1},a_j)(m,1)(a_j,a_{j-1}).
\end{split}
\end{equation}

The second of equalities (\ref{E6}) has an analogous  form:
\begin{equation}\label{E8}
\begin{split}
&(m-2)(a_i,a_{i+1})^{p_0-1}(m-1)(a_i,a_{i+1})\dots
(m-2)(a_{j-1},a_j)^{p_{i-j+1}-1}\cdot \\
&m(a_{j-1},a_j)= m(a_{i+1},a_i)(m,2)(a_{i+1},a_i)^{p_0-1}\dots
(m,1)(a_j,a_{j-1})\cdot \\
&(m,2)(a_j,a_{j-1})^{p_{i-j+1}-1}.
\end{split}
\end{equation}

Also we will use the following Garside's formula:
\begin{equation}\label{E9}
a_sa_{i+1j}=a_{i+1j}a_{s+1}, j\leq i,j\leq s\leq i.
\end{equation}
This formula looks like
\begin{align*}
&m(a_s,a_{i+1})\dots(m-1)(a_s,a_{s+1})\dots(m-1)(a_{s+1},a_j)=\\
&(m,1)(a_{i+1},a_s)\dots(m,1)(a_{s+1},a_s)\dots m(a_j,a_{s+1}).
\end{align*}

By $W(a_j,\dots,a_i)=W(j,i)$ we will mean any word in alphabet
${a_j,a_{j+1},\dots,a_i}$ if $j\leq i$, and empty word if $j>i$.

\begin{thm}\label{T1}
A \GS\ bases of $B_{n+1}^+$ consists of the following
relations :
\begin{align}
&a_{i+1}a_iV(1,i-1)W(j,i)a_{i+1j}=a_ia_{i+1}a_iV(1,i-1)a_{ij}W'(j+1,i+1),\label{E10}\\
 &a_sa_k=a_ka_s, s-k\geq 2,\label{E11}
\end{align}
where $1\leq i\leq n-1,$ $1\leq j\leq i+1, ,$ $
W'=W(a_{j+1},\dots, a_{i+1})$, and $W$ begins with $a_i$ if it is
not empty.
\end{thm}

\begin{rem}
In    particular the following relations are in the previous list:
\begin{align*}
&a_{i+1}a_ia_{i+1}=a_ia_{i+1}a_i\mbox{  } (V=1,j=i+1),\\
 &a_{i+1}a_i^la_{i+1}a_i=a_ia_{i+1}a_i^2a_{i+1}^{l-1},l\geq 2\mbox{   }
(V=1,W=a_i^{l-1},j=i),\\
&a_{i+1}a_i^la_{i-1j}(\mathbf{p})a_{i+1j}=a_ia_{i+1}a_i^2a_{i-1j}a_{i+1}^{l-1}a_{ij+1}(
\mathbf{p}),l\geq 2,i-1\geq j \\
&(V=1,W=a_i^{l-1}a_{i-1j}(\mathbf{p})).\\
\end{align*}

This relations have a form (\ref{E3}) because of Lemma \ref{L1}
and (\ref{E8}). In general (if $W$ is not empty) relation
(\ref{E10}) has a form (\ref{E4}). It follows from above and
(\ref{E9}).

Also the rite hand sides of (\ref{E6}) go to the left hand sides
using the eliminations of leading words of relations (\ref{E10}).
\end{rem}
\begin{proof}

Let $\mathbf{S}$ is the set of relations (polynomials) from the
theorem. We need to proof that any composition of polynomials from
$\mathbf{S}$ is trivial $mod \mathbf{S}.$

First of all let us deal with compositions of including. It is
easy to see that any composition of including of (\ref{E11}) in
(\ref{E10}) is trivial for the only possible composition of
this including should correspond of some including of
$a_sa_k,s-k\geq 2$ in $VW$, or else of $a_ia_k,i-k\geq 2$ in
$a_iV$. So we need to figure out any composition of including
relations of the type (\ref{E10}). Let $f$ be a relation
(polynomial) (\ref{E10}) and $g$ be the following relation:
$$
a_{i_1+1}a_{i_1}V_1(1,i_1-1)W_1(j_1.i_1)a_{i_1+1j_1}=
a_{i_1}a_{i_1+1}a_{i_1}V_1(1,i_1-1)a_{i_1j_1}W_1'(j_1+1.i_1+1),
$$
where $j_1\leq i_1+1$ and $W_1$ begins with $a_{i_1}$ if not
empty. Let $u,u_1$ be the leading monomials of $f,g$ respectively.
Suppose that $u_1$ is a subword of $u$. Following possibilities
may occur:

1)$i=i_1$. Then $V=V_1,W=W_1$, $j_1\geq j.$

Then
\begin{align*}
&(f,g)_u=a_ia_{i+1}a_iVa_{ij_1}W'a_{j_1-1j}-a_ia_{i+1}a_iVa_{ij}W'\equiv0,\\
\end{align*}
for $W'=W_1'=W(a_{j_1+1},\dots,a_{i+1})$, hence
$W'a_{j_1-1j}\equiv a_{j_1-1j}W'$.

2)$i=i_1+1.$ Then
$$
u_1=a_ia_{i-1}V_1(1,i-2)W_1(j_1,i-1)a_{ij_1}.
$$
It means that
$$
V(1,i-1)=a_{i-1}V_1(1,i-2)W_1(j_1,i-1),W(j,i)=a_{ij_1}W_2(j,i).
$$
We have
\begin{align*}
&(f,g)_u=a_{i+1}a_{i-1}a_ia_{i-1}V_1(1,i-2)a_{i-1j_1}W_1'(j_1+1,i)W_2(j,i)a_{i+1j}-\\
&a_ia_{i+1}[a_ia_{i-1}V_1(1,i-2)W_1(j_1,i-1)a_{ij_1}]a_{j_1-1j}a_{i+1j_1+1}
W_2'(j+1,i+1)\equiv\\
&a_{i-1}[a_{i+1}a_ia_{i-1}V_1(1,i-2)a_{i-1j_1}W_1'(j_1+1,i)W_2(j,i)a_{i+1
j}]-\\
&a_ia_{i+1}a_{i-1}a_ia_{i-1}V_1(1,i-2)a_{i-1j_1}W_1'(j_1+1,i)a_{j_1-1j}a_
{i+1j_1+1}
W_2'(j+1,i+1)\equiv\\
&a_{i-1}a_ia_{i+1}[a_ia_{i-1}V_1(1,i-2)a_{i-1j_1}a_{ij_1}]a_{j_1-1j}W_1''
(j_1+2,i+1)
W_2'(j+1,i+1)-\\
&a_ia_{i-1}[a_{i+1}a_ia_{i-1}V_1(1,i-2)a_{i-1j_1}a_{j_1-1j}W_1'(j_1+1,i)a
_{i+1j_1+1}]
W_2'(j+1,i+1)\equiv\\
&a_{i-1}a_ia_{i+1}a_{i-1}a_ia_{i-1}V_1(1,i-2)a_{i-1j_1}a_{ij_1+1}a_{j_1-1
j}W_1''(j_1+2,i+1)
W_2'(j+1,i+1)-\\
&a_ia_{i-1}a_ia_{i+1}a_ia_{i-1}V_1(1,i-2)a_{i-1j}a_{ij_1+1}W_1''(j_1+2,i+
1)
W_2'(j+1,i+1)\equiv 0,
\end{align*}
for $a_{ij_1+1}a_{j_1-1j}\equiv a_{j_1-1j}a_{j_1-1j}$.

3) $i-i_1\geq 2.$ The $u_1$ is a subword of
$V(1,i-1)W(j,i)a_{i+1j}$. Than it should be a subword either $V$
or $W$. First case is clear. Let us figure out second case:
\begin{align*}
&W=W_2(j,i)a_{i_1+1}a_{i_1}V_1(1,i_1-1)W_1(j_1,i_1)a_{i_1+1j_1}W_3(j,i)
,\\
&w=a_{i+1}a_iV(1,i-1)W_2(j,i)a_{i_1+1}a_{i_1}V_1(1,i_1-1)W_1(j_1,i_1)
a_{i_1+1j_1}W_3(j,i)a_{i+1j},\\
&(f,g)_w=a_ia_{i+1}a_iV(1,i-1)a_{ij}W_2'(j+1,i+1)a_{i_1+2}a_{i_1+1}V_1'
(2,i_1)
W_1'(j_1+1,i_1+1)\cdot\\
&a_{i_1+2j_1+1}W_3'(j+1,i+1)-[a_{i+1}a_iV(1,i-1)W_2(j,i)a_{i_1}a_{i_1+1}a
_{i_1}
V_1(1,i_1-1)a_{i_1+1j_1}\cdot\\
&W_1'(j_1+1,i_1+1)W_3(j,i)a_{i+1j}]\equiv
a_ia_{i+1}a_iV(1,i-1)a_{ij}W_2'(j+1,i+1)[a_{i_1+2}a_{i_1+1}\cdot\\
&V_1'(2,i_1)W_1'(j_1+1,i_1+1)a_{i_1+2j_1+1}]W_3'(j+1,i+1)-a_ia_{i+1}a_iV(
1,i-1)
a_{ij}\cdot\\
&W_2'(j+1,i+1)a_{i_1+1}a_{i_1+2}a_{i_1+1}V_1'(2,i_1)a_{i_1+1j_1+1}
W_1''(j_1+2,i_1+2)W_3'(j+1,i+1) \\
&\equiv0.
\end{align*}

Now we need to check compositions of intersection.

 It is easy to see that all compositions of
intersection with (\ref{E11}) are trivial. Therefore we only need
to consider compositions of intersection of (\ref{E10}) relative
to the following ambiguities $w$:

1)It is a case of one letter intersection,
\begin{align*}
&
w=a_{i+1}a_iV(1,i-1)W(j,i)a_{i+1j}a_{j-1}V_1(1,j-2)W_1(j_1,j-1)a_{jj_1}
,\\
&(f,g)_w=a_ia_{i+1}a_iV(1,i-1)a_{ij}W'(j+1,i+1)a_{j-1}V_1(1,j-2)
W_1(j_1,j-1)a_{jj_1}-\\
&a_{i+1}a_iV(1,i-1)W(j,i)a_{i+1j+1}a_{j-1}a_ja_{j-1}V_1(1,j-2)a_{j-1j_1}
W_1'(j_1+1,j)\equiv \\
&a_ia_{i+1}a_iV(1,i-1)a_{ij-1}W'(j+1,i+1)V_1(1,j-2)
W_1(j_1,j-1)a_{jj_1}-\\
&[a_{i+1}a_iV(1,i-1)(W(j,i)a_{j-1})a_{i+1j-1}]V_1(1,j-2)a_{j-1j_1}
W_1'(j_1+1,j)\equiv \\
&a_ia_{i+1}a_iV(1,i-1)a_{ij-1}W'(j+1,i+1)V_1(1,j-2)
W_1(j_1,j-1)a_{jj_1}-\\
&a_ia_{i+1}a_iV(1,i-1)a_{ij-1}(W'(j+1,i+1)a_j)V_1(1,j-2)a_{j-1j_1}
W_1'(j_1+1,j)\equiv0,
\end{align*}
for $a_jV_1(1,j-2)\equiv a_jV_1(1,j-2)$ and
$W_1(j_1,j-1)a_{jj_1}\equiv a_{jj_1}W_1'(j_1+1,j)$;

 2) It is a case of two letters intersection,
\begin{align*}
&
w=a_{i+1}a_iV(1,i-1)W(j,i)a_{i+1j+2}a_{j+1}a_jV_1(1,j-1)W_1(j_1,j)a_{j+
1j_1},\\
&(f,g)_w=a_ia_{i+1}a_iV(1,i-1)a_{ij}W'(j+1,i+1)V_1(1,j-1)
W_1(j_1,j)a_{j+1j_1}-\\
&a_{i+1}a_iV(1,i-1)W(j,i)a_{i+1j+2}a_ja_{j+1}a_jV_1(1,j-2)a_{jj_1}
W_1'(j_1+1,j+1)\equiv \\
&a_ia_{i+1}a_iV(1,i-1)a_{ij}W'(j+1,i+1)V_1(1,j-1)
W_1(j_1,j)a_{j+1j_1}-\\
&[a_{i+1}a_iV(1,i-1)(W(j,i)a_j)a_{i+1j}]V_1(1,j-2)a_{jj_1}
W_1'(j_1+1,j+1)\equiv \\
&a_ia_{i+1}a_iV(1,i-1)a_{ij}W'(j+1,i+1)V_1(1,j-1)
W_1(j_1,j)a_{j+1j_1}-\\
&a_ia_{i+1}a_iV(1,i-1)a_{ij}(W'(j+1,i+1)a_{j+1})V_1(1,j-2)a_{jj_1}
W_1'(j_1+1,j+1)\equiv 0,\\
\end{align*}
for $a_{j+1}V_1(1,j-2)\equiv V_1(1,j-2)a_{j+1}$ and
$W_1(j_1,j)a_{j+1j_1}\equiv a_{j+1j_1}W_1'(j_1+1,j+1).$

3) It is a case of more than two letters intersection:

\begin{align*}
& w=a_{i+1}a_iV(1,i-1)W(j,i)a_{i+1}a_i\dots a_{i_1+1}a_{i_1}\dots
a_jV_2(1,i_1-1)W_1(j_1,i_1)a_{i_1+1j_1},\\
&(f,g)_w=a_ia_{i+1}a_iV(1,i-1)a_{ij}W'(j+1,i+1)V_2(1,i_1-1)W_1(j_1,i_1)=
a_{i_1+1j_1}-\\
&a_{i+1}a_iV(1,i-1)W(j,i)a_{i+1i_1+2}a_{i_1}a_{i_1+1}a_{i_1}a_{i_1-1j}
V_2(1,i_1-1)a_{i_1j_1}
W_1'(j_1+1,i_1+1)\equiv \\
&a_ia_{i+1}a_iV(1,i-1)a_{ij}W'(j+1,i+1)V_2(1,i_1-1)W_1(j_1,i_1)a_{i_1+1j_=
1}-\\
&[a_{i+1}a_iV(1,i-1)(W(j,i)a_{i_1})a_{i+1j}]V_2(1,i_1-1)a_{i_1j_1}
W_1'(j_1+1,i_1+1)\equiv \\
&a_ia_{i+1}a_iV(1,i-1)a_{ij}W'(j+1,i+1)V_2(1,i_1-1)W_1(j_1,i_1)a_{i_1+1j_=
1}-\\
&a_ia_{i+1}a_iV(1,i-1)a_{ij}(W'(j+1,i+1)a_{i_1+1})V_2(1,i_1-1)a_{i_1j_1}
W_1'(j_1+1,i_1+1)\equiv 0,
\end{align*}
for $a_{i_1+1}V_2(1,i_1-1)\equiv V_2(1,i_1-1)a_{i_1+1}$ and
$W_1(j_1,i_1)a_{i_1+1j_1}\equiv a_{i_1+1j_1}W_1'(j_1+1,i_1+1).$

We have checked all compositions and the theorem is proved.

\end{proof}

\begin{rem}
In  view of \cite{Ga69}, it is important to have an efficient
algorithm to know whether a word $v$ of $B_{n+1}$ begin with the
fundamental word $ \Delta$ in some presentation, where
$$
\Delta = \Lambda_1 \Lambda_2 \dots \Lambda_n, \Lambda_i=a_i\dots
a_1, i\geq 1.
$$

Let us formulate the following statement that easily follows from
the previous theorem.

\begin{prop}\label{T2}
Let $v\in PBW(B_{n+1}^+)$. Then $v$ begins with $\Delta$ in
some presentation if and only if
$$
v\doteq \Lambda_1 V_1(1) \Lambda_2 V_2(1,2)\dots
\Lambda_{n-1}V_{n-1}(1,n-1)\Lambda_nu,
$$
where $V_i(1,i)$ begins with $a_1 $ or empty. In this case we have
$$
v=\Delta V_1^{(n-1)}(n)V_2^{(n-2)}(n-1,n)\dots V_{n-1}'(2,n).
$$
\end{prop}
\begin{proof}The result follows imediatly from the following formula:
$$
PBW(a_{i+1}a_i\dots a_1w)\doteq V(1,i)a_{i+1}a_i\dots a_1w', w,w'\in B_{n+1}^+.
$$
\end{proof}
\end{rem}
\begin{rem}
In  view of \cite{Ga69}, Theorem \ref{T1} gives rise to a new algorithm of the solution
of the word problem for the braid group $B_{n+1}$. Indeed, any group word $w\in B_{n+1} $
can be presented in a form $ \Delta ^kv  $, where $k\in \mathbb{Z}$,  $ v\in  B_{n+1}^+$.
Using Proposition
\ref{T2}, we may present $v$ in a form
$$
v=\Delta^su,
$$
where $s\in \mathbb{N} $, $u$ does not begin with $\Delta$ in any presentation
(it means that $u\neq \Delta u_1$ in the braid semigroup). Also we may assume
that $u$ is a reduced word, $u\in PBW(  B_{n+1}^+)$.
Now $u$ is nothing else as $ \overline{v} $ in the sense
of \cite{Ga69}, Theorem 5. As a result, the presentation
$$
w=\Delta^{k+s}u
$$
is unique due to the same Theorem 5 \cite{Ga69}. Let us call the last
presentation the Garsaid normal form of an element of the braid group.
Garside normal form is effective due to Theorem \ref{T1}. Now the new
algorithm for the solution of the word problem for the braid group
can be written down using the ELW's of the \GS\ basis of this group.
\end{rem}

\end{document}